\documentclass[journal]{IEEEtran}
\ifCLASSINFOpdf
\else
\fi
\usepackage{amsmath}
\usepackage{bm}
\usepackage{calc}
\usepackage{esint}
\usepackage{graphicx}
\usepackage{multicol}
\usepackage{multirow}
\usepackage{tikz}
\usepackage{xcolor}
\usepackage{cite}

\usepackage{soul}
\usepackage[flushleft]{threeparttable}
\usepackage{upgreek}
\usepackage[nice]{nicefrac}
\makeatletter
\newcommand*\bigcdot{\mathpalette\bigcdot@{.5}}
\newcommand*\bigcdot@[2]{\mathbin{\vcenter{\hbox{\scalebox{#2}{$\m@th#1\bullet$}}}}}
\makeatother

\newcommand{\dt}{\Delta t}
\newcommand{\dx}{\Delta x}
\newcommand{\dy}{\Delta y}
\newcommand{\dz}{\Delta z}

\definecolor{darkRed}{rgb}{0.7, 0, 0}
\definecolor{darkBlue}{rgb}{0, 0, 0.7}

\begin{document}
%
\title{A Dissipation Theory for Potentials-Based FDTD for Lossless Inhomogeneous Media}
%
%
%

\author{Fadime~Bekmambetova,~\IEEEmembership{Student Member,~IEEE}
        and Piero~Triverio,~\IEEEmembership{Senior Member,~IEEE}
\thanks{Manuscript received ...; revised ....}
\thanks{This work was supported in part by the Natural
	Sciences and Engineering Research Council of Canada (Discovery Grant Program) and in part by the Canada Research Chairs Program.}
\thanks{F. Bekmambetova and P. Triverio are with the Edward S. Rogers Sr. Department of Electrical \& Computer Engineering, University of Toronto, Toronto, ON, M5S 3G4 Canada, e-mails: fadime.bekmambetova@mail.utoronto.ca,
	piero.triverio@utoronto.ca. (Corresponding author: Piero Triverio.)}
\thanks{\copyright 2021 IEEE. Personal use of this material is permitted. Permission from IEEE must be obtained for all other uses, in any current or future media, including reprinting/republishing this material for advertising or promotional purposes, creating new collective works, for resale or redistribution to servers or lists, or reuse of any copyrighted component of this work in other works.}
}

%
%

\markboth{IEEE Antennas and Wireless Propagation Letters}%
{Shell \MakeLowercase{\textit{et al.}}: Bare Demo of IEEEtran.cls for IEEE Journals}
%




\maketitle

\begin{abstract}
A dissipation theory is proposed for the potentials-based FDTD algorithm for the case of inhomogeneous lossless media. We show that under the Courant-Friedrichs-Lewy (CFL) limit, the equations describing the time evolution of scalar and vector potentials can be seen as a lossless system. The developed theory provides insights into how electromagnetic energy and power flow are approximated in FDTD schemes. It can also be used to create new algorithms with guaranteed stability.
\end{abstract}


%
\IEEEpeerreviewmaketitle

\section{Introduction}
%
%
%
%
\IEEEPARstart{T}{he} finite-difference time-domain (FDTD) algorithm~\cite{yee,gedney-2011} is widely used for solving Maxwell's equations. Recently, there has been an increased interest in integrating FDTD with quantum models~\cite{sui-2007-ims,ahmed-2010-awpl,ohnuki-2013-ijnm,takeuchi-2015-pr,ryu-fdtd-quantum-paper-2016,xiang-2019-ieee-access}. Interaction between electromagnetic signals and quantum states plays an essential role in many architectures being considered for quantum computing~\cite{jazaeri-qc-review}, calling for algorithms that can jointly simulate quantum and electromagnetic phenomena. 
Quantum models of field-particle interaction typically involve potentials, as opposed to the fields computed in the traditional FDTD. The requirement of the knowledge of potentials in quantum modeling makes potentials a natural choice of unknowns for FDTD in quantum applications~\cite{ryu-fdtd-quantum-paper-2016,xiang-2019-ieee-access}. Earlier, potentials-based FDTD~(P-FDTD) formulations have been investigated, for example, as means of reducing computational requirements~\cite{kunz-fdtd-for-em-1993,flaviis-1998-mgwl}.


P-FDTD methods still lack many of the advancements that have been proposed for the fields-based FDTD, including  subgridding~\cite{okoniewski-3d-subgridding-1997}, model order reduction~\cite{kulas-2004-mwcl}, and many others. One of the difficulties in creating such new schemes is the complexity of ensuring stability.
In the case of the traditional fields-based FDTD, one needs to select the time step below the Courant-Friedrichs-Lewy (CFL) limit~\cite{kim-von-neumann,gedney-2011}
\begin{equation}
	\label{eq:cfl}
	\dt <  \dfrac{\sqrt{\mu\varepsilon}}{\sqrt{(\dx)^{-2}+(\dy)^{-2}+(\dz)^{-2}}}
	\,.
\end{equation}
The same stability requirement has been used for P-FDTD algorithms in~\cite{ryu-fdtd-quantum-paper-2016,flaviis-1998-mgwl,xie-2021-tps}. 
Various techniques have been proposed for ensuring stability of new FDTD schemes, such as von Neumann analysis~\cite{kim-von-neumann}, the iteration matrix method~\cite{remis-2000-jcp}, as well as methods based on analyzing a linear operator involved in the second order formulation of FDTD equations~\cite{mrozowski-1994-mgwl,kulas-2004-mwcl}. 

In this work, we propose a dissipation-based framework for ensuring stability of P-FDTD schemes, extending previous works on traditional FDTD~\cite{jnl-2018-fdtd-3d-dissipative}. The dissipation framework takes root in the theory of dissipative systems~\cite{willems1972dissipative} and is a generalization of the energy method~\cite{edelvik2004general}. 
The present work provides proper expressions to compute the energy stored in a region and power absorbed through its boundary in P-FDTD. It also gives insights into how electromagnetic energy is approximated in P-FDTD.	
The proposed approach can facilitate the stability analysis in complex setups involving many parts, such as grids of different resolution in subgridding~\cite{okoniewski-3d-subgridding-1997}, different types of equations in hybrid methods~\cite{bretones1998hybrid} or in multi-physics simulations~\cite{meng-2012-jctc,ryu-fdtd-quantum-paper-2016}.

\section{P-FDTD Equations in State Space Form}
\label{sec:state-space-form}

We consider a rectangular region filled with an inhomogeneous isotropic dielectric with permittivity $\varepsilon$ and permeability $\mu$. In this section, we discretize the wave equations governing the scalar potential $\phi$ and vector potential $\vec A$ and cast the discretized equations in the form of a state space dynamical system suitable for analyzing energy dissipation.

\subsection{The Continuous State Equations for $\phi$ and $\vec A$}
Under the generalized Lorenz gauge~\cite{nisbet-gen-lorenz-gauge-1957,flaviis-1998-mgwl,chew-gen-gauge-pier}, the wave equation for scalar potential is
\begin{equation}
\label{eq:phi-second-order}
\chi \dfrac{\partial^2 \phi}{\partial t^2} 
= \nabla \bigcdot \varepsilon \nabla \phi
\end{equation}
where $\chi = \mu \varepsilon^2$. Equation~\eqref{eq:phi-second-order}
can be written in a form involving only first order time derivatives by introducing the following variables: $[\nabla \phi]$ for the gradient of scalar potential and $[\partial_t \phi]$ for the temporal derivative of scalar potential. The square brackets are used to indicate that we take the values of $\nabla \phi$ and  $\partial_t \phi$ as unknowns. Hence, \eqref{eq:phi-second-order} can be rewritten as
\begin{subequations}
	\label{eq:phi-state-space}
	\begin{eqnarray}	
	\label{eq:der-phi-state-space}
	\chi \dfrac{\partial}{\partial t}[\partial_t \phi] &=& \nabla \bigcdot \varepsilon [\nabla \phi]
	\\
	\label{eq:grad-phi-state-space}
	\varepsilon \dfrac{\partial}{\partial t} [\nabla \phi] &=& \varepsilon \nabla [\partial_t \phi] \,.
	\end{eqnarray}
\end{subequations}
Similar equations can be derived to describe the temporal evolution of $\vec A$
\begin{subequations}
	\label{eq:a-state-space}
	\begin{eqnarray}
	\label{eq:a-state-space-b}
	\mu^{-1} \dfrac{\partial \vec B}{\partial t} &=& \mu^{-1} \nabla \times [\partial_t \vec A]
	\\
	\label{eq:a-state-space-der-a}
	\varepsilon \dfrac{\partial [\partial_t \vec A]}{\partial t} &= &-\nabla \times \mu^{\!-1} \vec B -\varepsilon \nabla \kappa
	\\
	\label{eq:a-state-space-kappa}
	\chi \dfrac{\partial \kappa}{\partial t} &=& -\nabla \bigcdot \varepsilon [\partial_t \vec A]
	\end{eqnarray}
\end{subequations}
where $\vec B$ is magnetic flux density and $\kappa$ is defined as $-\chi^{-1} \nabla \bigcdot \varepsilon \vec A$. Equations~\eqref{eq:phi-state-space} and \eqref{eq:a-state-space} are equivalent to the wave equation for $\vec A$ and $\phi$ involving second order time derivatives~\cite{chew-gen-gauge-pier}.

\subsection{Discretized State Equations}
\label{sec:discretization}

\begin{figure}[t!]
	\centering
	\includegraphics[scale=1]{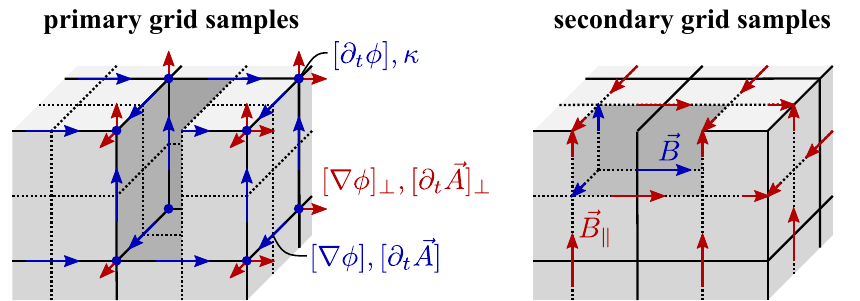}
	\caption{Location of conventional P-FDTD variables (blue) and hanging variables (red). The dotted lines show the secondary grid. The cutouts show where samples are located inside the region.}
	\label{fig:region}
\end{figure}

\begin{figure}[t!]
	\centering
	\includegraphics[scale=1]{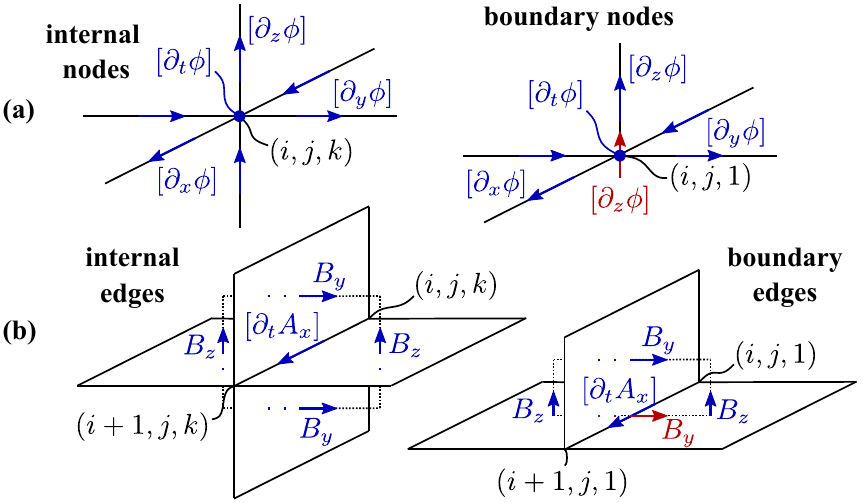}
	\caption{(a) Samples involved in equations \eqref{eq:der-phi-update-inner} and \eqref{eq:der-phi-update-k=1}. (b) Samples involved in equations discretizing \eqref{eq:a-state-space-der-a}, such as~\eqref{eq:der-a-update-k=1}.}
	\label{fig:sampling}
\end{figure}

We divide the region into $n_x \times n_y \times n_z$ primary cells with dimension $\dx \times \dy \times \dz$ and discretize equations \eqref{eq:phi-state-space} and \eqref{eq:a-state-space} using the centered finite differences. Overall, the discretization follows the approach in \cite{ryu-fdtd-quantum-paper-2016}, with differences in the treatment of the region boundary, as will be explained in the subsequent discussion. The variables are sampled on the staggered primary and secondary grids~\cite{gedney-2011}, as illustrated in Fig.~\ref{fig:region}. The secondary cells strictly inside the region have dimension $\dx \times \dy\times \dz$, whereas the cells adjacent to the boundary are halved in one of the dimensions. For cells adjacent to two or three boundary faces, two and three of the dimensions are halved, respectively.
The samples of the variables in~\eqref{eq:phi-state-space}--\eqref{eq:a-state-space} located inside the region are updated with equations equivalent to~\cite{ryu-fdtd-quantum-paper-2016}. However, on the region boundary, special equations must be used to create a self-contained model that does not involve samples outside the region, which allows precisely quantifying the amount of electromagnetic energy stored in the region.

\subsubsection{Scalar Potential}
The variables $[\nabla \phi]$ and $[\partial_t \phi]$ are sampled at the edges and nodes of the primary grid, respectively, as shown in Fig.~\ref{fig:region} in blue. By approximating the gradient operator in \eqref{eq:grad-phi-state-space} via centered finite differences, we obtain
\begin{multline}
	\label{eq:grad-phi-update-inner}
\dx \dy \dz \dfrac{\varepsilon}{\dt} \left(
[\partial_x\phi]_{i+\frac12, j,k}^{n+1}
- [\partial_x\phi]_{i+\frac12, j,k}^n
\right) = \\
\varepsilon \dy \dz \left(
[\partial_t\phi]_{i+1,j,k}^{n+\frac12} 
- [\partial_t\phi]_{i,j,k}^{n+\frac12}
\right)
\end{multline}
for the $x$-directed edges strictly inside the region.  Equation~\eqref{eq:grad-phi-update-inner} can be seen as a discrete integration of \eqref{eq:grad-phi-state-space} over the volume associated with the primary edge where the $[\partial_x\phi]$ sample is located. In particular, $\dx$ is the length associated with that primary edge and $\dy \dz$ is the area of the secondary face pierced by that edge. For edges on the bottom face of the boundary ($k=1$), we take the same equation as \eqref{eq:grad-phi-update-inner}, but with $\dz$ replaced by $\dz/2$, since the area of the corresponding secondary face is $\dy \dz/2$. 
Analogous modifications are made for the edges shared between two faces of the boundary.

With reference to Fig.~\ref{fig:sampling}a, \eqref{eq:der-phi-state-space} is discretized as
\begin{multline}
\label{eq:der-phi-update-inner}
\dx \dy \dz \dfrac{\chi}{\dt} \left(
[\partial_t\phi]_{i,j,k}^{n+\frac12}
-[\partial_t\phi]_{i,j,k}^{n-\frac12}
\right) \\
= \dy \dz \left(\varepsilon [\partial_x\phi]_{i+\frac12,j,k}^n - \varepsilon [\partial_x\phi]_{i-\frac12,j,k}^n\right)\\
+ \dx \dz \left(\varepsilon [\partial_y\phi]_{i,j+\frac12,k}^n - \varepsilon [\partial_y\phi]_{i,j-\frac12,k}^n\right) \\
+ \dx \dy \left(\varepsilon [\partial_z\phi]_{i,j,k+\frac12}^n - \varepsilon [\partial_z\phi]_{i,j,k-\frac12}^n\right)
\end{multline}
for the nodes strictly inside the region. On the right hand side of \eqref{eq:der-phi-update-inner} is the discrete flux of $\varepsilon \nabla \phi$ through the faces of the secondary cell that has the node $(i,j,k)$ at its center. 
For nodes on the $k=1$ boundary, \eqref{eq:der-phi-update-inner} would involve samples outside the boundary, which the model aims to avoid. Instead, we introduce a hanging variable $[\partial_z\phi]|_{i,j,1}$ on the boundary of the region~\cite{venkatarayalu2007stable,jnl-2018-fdtd-3d-dissipative}, as shown in Fig.~\ref{fig:sampling}a in red. The equation is written over the $\dx \times \dy \times \dz/2$ secondary cell adjacent to the $k=1$ boundary
\begin{multline}
\label{eq:der-phi-update-k=1}
\dx \dy \dfrac{\dz}2 \dfrac{\chi}{\dt} \left(
[\partial_t \phi]_{i,j,1}^{n+\frac12}
-[\partial_t \phi]_{i,j,1}^{n-\frac12}
\right) =\\
 \dy \dfrac{\dz}2 \left(\varepsilon [\partial_x\phi]_{i+\frac12,j,1}^n - \varepsilon [\partial_x\phi]_{i-\frac12,j,1}^n\right)\\
+ \dx \dfrac{\dz}2 \left(\varepsilon [\partial_y\phi]_{i,j+\frac12,1}^n - \varepsilon [\partial_y\phi]_{i,j-\frac12,1}^n\right) \\
+ \dx \dy \varepsilon [\partial_z\phi]_{i,j,\frac32}^n 
- \dx \dy \varepsilon [\partial_z\phi]_{i,j,1}^n\,.
\end{multline}
Equations for the nodes shared between two or three sides of the boundary involve, respectively, two and three hanging variables. The hanging variables will be crucial to quantify the power absorbed by the region through the boundary~\cite{jnl-2018-fdtd-3d-dissipative}.

\subsubsection{Vector potential}
\label{sec:a-discretization}

As shown in Fig.~\ref{fig:region}, $\vec B$ is sampled on the secondary edges. The variables $[\partial_t \vec A]$ and $\kappa$ are sampled on the primary edges and primary nodes, respectively.
Equations~\eqref{eq:a-state-space-b}--\eqref{eq:a-state-space-kappa} are discretized analogously to \eqref{eq:der-phi-state-space}--\eqref{eq:grad-phi-state-space}. For edges and nodes strictly inside the region, the discrete equations are equivalent to those in \cite{ryu-fdtd-quantum-paper-2016}. In order to discretize equation \eqref{eq:a-state-space-kappa} on boundary nodes, we introduce a hanging variable for the normal component of the vector potential. For example, the equation on the $k=1$ boundary is written in an analogous way to~\eqref{eq:der-phi-update-k=1}, involving a hanging variable $[\partial_tA_z]|_{i,j,1}$.	
Likewise, the discretizated equation \eqref{eq:a-state-space-der-a} for $[\partial_t \vec A]$ on primary edges tangential to the boundary requires a hanging variable for the tangential magnetic flux density. For instance, on the $k=1$ boundary, the resulting equation reads
\begin{multline}
\label{eq:der-a-update-k=1}
\dx \dy \dfrac{\dz}2 \dfrac{\varepsilon}{\dt} \left(
	[\partial_t A_x]_{i+\frac12,j,1}^{n+1} 
	- [\partial_t A_x]_{i+\frac12,j,1}^{n}
	\right) = \\
- \dx \bigg(
\dfrac{\dz}2 \mu^{-1} B_z|_{i+\frac12,j+\frac12, 1}^{n+\frac12}
- \dfrac{\dz}2 \mu^{-1} B_z|_{i+\frac12,j-\frac12, 1}^{n+\frac12}
\\
-\dy \mu^{-1} B_y|_{i+\frac12,j,\frac32}^{n+\frac12}
+ \dy \mu^{-1} B_{y}|_{i+\frac12,j,1}^{n+\frac12}
\bigg) \\
- \dy \dfrac{\dz}2 \varepsilon \left(
\kappa|_{i+1,j,1}^{n+\frac12}
- \kappa|_{i,j,1}^{n+\frac12}
\right)
\end{multline}
The expression in the first bracket on the right hand side approximates the circulation of the magnetic field around the path in Fig.~\ref{fig:sampling}b.
The hanging variable for tangential boundary magnetic flux density $B_{y}|_{i+0.5,j,1}$ is introduced to avoid using samples outside of the region and later devise an expression for the energy entering the region through the boundary.

\subsection{State Equations in Matrix Form}
\newcommand{\volSec}{\mathbf{\Lambda}_V''}
\newcommand{\CHI}{\mathbf{\Lambda}_{\chi}}
\newcommand{\EPS}{\mathbf{\Lambda}_{\varepsilon}}
\newcommand{\MUinv}{\mathbf{\Lambda}_{\nicefrac1{\mu}}}
\newcommand{\derPhi}{[\partial_t \bm{\upphi}]}
\newcommand{\gradPhi}{[\nabla \bm{\upphi}]}
\newcommand{\D}{\mathbf{D}}
\newcommand{\C}{\mathbf{C}}
\newcommand{\areaSec}{\mathbf{\Lambda}_{S}''}
\newcommand{\areaPrim}{\mathbf{\Lambda}_{S}'}
\newcommand{\lenPrim}{\mathbf{\Lambda}_{l}'}
\newcommand{\lenSec}{\mathbf{\Lambda}_{l}'}

\newcommand{\Dnor}{\D_{n}}
\newcommand{\areaSecBoundary}{\mathbf{\Lambda}_{S,b}''}
\newcommand{\EPSnor}{\mathbf{\Lambda}_{\varepsilon,\bot}}
\newcommand{\gradPhiNor}{\gradPhi_{\bot}}
\newcommand{\derPhiBoundary}{[\partial_t \bm{\upphi}]_b}

\newcommand{\Ctan}{\mathbf{C}_t}
\newcommand{\lenSecTan}{\mathbf{\Lambda}_{l,t}'}
\newcommand{\MUinvTan}{\mathbf{\Lambda}_{\nicefrac1{\mu},t}}
\newcommand{\bbTan}{\mathbf b_t}

\newcommand{\bb}{\mathbf b}
\newcommand{\Kappa}{\bm {\upkappa}}
\newcommand{\derA}{[\partial_t \bm{a}]}

\newcommand{\KappaBoundary}{\bm {\upkappa}_b}
\newcommand{\derAtan}{[\partial_t \bm{a}]_t}
\newcommand{\derAnor}{[\partial_t \bm{a}]_n}

\newcommand{\Lnor}{\mathbf L_{\bot}}

\newcommand{\SignNor}{\mathbf{\Lambda}_{(\hat n \bigcdot)}}
\newcommand{\SignTan}{\mathbf{\Lambda}_{(\hat n \times)}}

State equations for scalar potential, such as \eqref{eq:grad-phi-update-inner}--\eqref{eq:der-phi-update-k=1}, written out for all internal and boundary nodes and edges can be collected in convenient matrix form
\begin{subequations}
\begin{equation}
\label{eq:compact-grad-phi}
\areaSec \lenPrim \dfrac{\EPS}{\dt} \left(\gradPhi^{n+1}-\gradPhi^n\right)
= 
-\areaSec \EPS \D^T \derPhi^{n+\frac12}
\end{equation}	
\begin{multline}
\label{eq:compact-der-phi}
	\volSec \dfrac{\CHI}{\dt} \left(\derPhi^{n+\frac12} - \derPhi^{n-\frac12}\right)
	= \D \areaSec \EPS \gradPhi^{n} \\+ \mathbf L_{\bot} \mathbf{\Lambda}_{(\hat n \bigcdot)} \areaSecBoundary \EPSnor \gradPhiNor^{n}
\end{multline}
\end{subequations}
where $\derPhi$ is a column vector collecting all $[\partial_t \phi]$ samples shown in blue in Fig.~\ref{fig:region} and Fig.~\ref{fig:sampling}a. Vector $\gradPhi$ collects the conventional $[\nabla \phi]$ samples, such as $[\partial_x \phi]$, also shown in blue in the figures. Vector $\gradPhiNor$ collects the corresponding hanging variables, shown in red in the figures. Symbol $\mathbf{\Lambda}$ denotes diagonal matrices. In particular, $\areaSec$ and $\lenPrim$ contain secondary face areas and primary edge lengths. Similarly, $\volSec$ contains volumes of secondary cells associated with the samples in $\derPhi$. Matrices $\EPS$ and $\CHI$ contain, respectively, the values of $\varepsilon$ and $\chi$ on the primary edges and nodes. Matrix $\D$ contains $1$'s, $-1$'s, and $0$'s needed for the computation of discrete divergence or flux. Elements of matrix $\mathbf L_{\bot}$ are $0$'s and $1$'s, with $1$ if the corresponding hanging variable in $\gradPhiNor$ is located at the node corresponding to the row of $\mathbf L_{\bot}$. Matrix $\mathbf{\Lambda}_{(\hat n \bigcdot)}$ contains $1$ or $-1$ on the diagonal, depending on the sign of the dot product between the \emph{outward} normal vector $\hat n$ to the boundary and the Cartesian unit vector associated with the additional edge where the hanging variable is located. Matrices $\areaSecBoundary$ and $\EPSnor$ contain boundary face areas and permittivities on the edges associated with the hanging variables in $\gradPhiNor$. In an analogous way, one can write the discretized equations for vector potential in matrix form.

Equations~\eqref{eq:compact-grad-phi}--\eqref{eq:compact-der-phi} can be rewritten in the generalized state space form
\begin{subequations}
	\label{eq:dyn-sys}
\begin{eqnarray}
	\label{eq:dyn-sys-a}
(\mathbf{R} + \mathbf F) \mathbf x^{n+\frac12} 
&=& (\mathbf{R} - \mathbf F) \mathbf x^{n-\frac12} 
+ \mathbf B \mathbf u^{n}\\\label{eq:dyn-sys-b}
\mathbf y^{n-\frac12} &=& \mathbf L^T \mathbf x^{n-\frac12}
\end{eqnarray}
\end{subequations}
where the state, input, and output vectors are given by
\begin{equation}
\mathbf x^{n-\frac12} \!=\! \begin{bmatrix}
\gradPhi^n\\ \derPhi^{n-\frac12}
\end{bmatrix}\!,
\,
\mathbf u^n = \gradPhiNor^n,
\,
\mathbf y^{n-\frac12} = \derPhiBoundary^{n-\frac12}
\end{equation}
and $\derPhiBoundary$ is a vector collecting the $[\partial_t \phi]$ samples located on the boundary. Matrices $\mathbf R$, $\mathbf F$, $\mathbf B$, and $\mathbf L$ are given by
\begin{equation}
\label{eq:mtx-def}
\mathbf R =\! \begin{bmatrix}
\mathbf R_{11}\!\! & \!\mathbf R_{21}^T \\
\mathbf R_{21}\!\! & \!\mathbf R_{22}
\end{bmatrix}\!\!,
\,\,
\mathbf F =\! \begin{bmatrix}
\mathbf 0\!\!\!\! & \!\!\mathbf R_{21}^T \\
-\mathbf R_{21}\!\!\!\! & \!\!\!\mathbf 0
\end{bmatrix}\!\!,
\,\,
\mathbf B = \mathbf L \mathbf S,
\,\,
\mathbf L =\! \begin{bmatrix}
\mathbf 0 \\ \Lnor\!
\end{bmatrix}
\end{equation}
where $\mathbf R_{11} = (\dt)^{-1}\EPS \areaSec \lenPrim$, $\mathbf R_{21} = \frac{1}2\D \areaSec \EPS$, $\mathbf R_{22} = (\dt)^{-1} \CHI \volSec$, $\mathbf S = \SignNor \areaSecBoundary \EPSnor$.
One can also formulate a dynamical system of the form~\eqref{eq:dyn-sys} for the vector potential. Writing the FDTD equations in this form enables the application of many theorems from systems theory, including those that investigate whether the system is lossless or not.

\section{Dissipativity Analysis}
\label{sec:dissipativity}
\subsection{Scalar Potential}

A dynamical system in the form \eqref{eq:dyn-sys} is said to be \emph{lossless} with supply rate $s$ if there is a storage function $\mathcal E$, satisfying
\begin{subequations}
	\label{eq:dissip-ineq-all}
\begin{equation}
\label{eq:storage-nonnegative}
\mathcal E(\mathbf x^{n-\frac12}) \ge 0,\,\, \mathcal E(\mathbf 0) = 0 \quad \forall n
\end{equation}
\begin{equation}
\label{eq:dissip-ineq}
\mathcal E(\mathbf x^{n+\frac12})-\mathcal E(\mathbf x^{n-\frac12}) = s(\mathbf u^n, \mathbf y^{n-\frac12}, \mathbf y^{n+\frac12}) \quad \forall n
\end{equation}
\end{subequations}
for all $\mathbf u^n$~\cite{byrnes1994losslessness}.
As storage function and supply rate, we propose to use the following expressions, based on~\cite{jnl-2018-fdtd-3d-dissipative,edelvik2004general}
\begin{equation}
\label{eq:storage}
\mathcal E(\mathbf x^{n-\frac12}) = \dfrac{\dt}{2} \left(\mathbf x^{n-\frac12}\right)^T \mathbf R \mathbf x^{n-\frac12}
\end{equation}
\begin{equation}
\label{eq:supply}
s(\mathbf u^n, \mathbf y^{n-\frac12}, \mathbf y^{n+\frac12}) = 
\dt \dfrac{\left(\mathbf{y}^{n-\frac12} + \mathbf{y}^{n+\frac12}\right)^{\!T}\!\!\!}2 \mathbf S \mathbf u^{n}\,.
\end{equation}
When \eqref{eq:storage} is expanded, it is revealed to be a discrete counterpart of the continuous expression for the stored energy found in literature \cite{chew-qem-new-look-2,puthoff-energy}
\begin{equation}
\label{eq:storage-cont-phi}
\mathcal E_{\text{cont}} = \iiint_V 
\left[
\frac{\varepsilon}{2} \nabla \phi \bigcdot \nabla \phi + \frac{\chi}{2} \left(\dfrac{\partial \phi}{\partial t}\right)^{\!2}
\right]
dV\,.
\end{equation}
Expression~\eqref{eq:supply} is a discrete approximation of the continuous energy entering the region during one time step~\cite{puthoff-energy}
\begin{equation}
	\label{eq:supply-cont-phi}
s_{\text{cont}} = \int_{\Delta t}\oiint_S \varepsilon \dfrac{\partial \phi}{\partial t} \nabla \phi \bigcdot \hat n\ dS\ dt
\end{equation}
where $-\varepsilon \frac{\partial \phi}{\partial t} \nabla \phi$ is the power flux density similar to the Poynting vector.

Using the structure of the matrices in~\eqref{eq:mtx-def}, we were able to prove that the system~\eqref{eq:dyn-sys} satisfies \eqref{eq:dissip-ineq-all} if $\mathbf R$ is positive definite.
For homogeneous media, it can be shown that ensuring the CFL limit \eqref{eq:cfl} guarantees the positive definiteness of $\mathbf R$ and thus ensures that \eqref{eq:dyn-sys} conserves energy.

Hence, in addition to the well-accepted connection to stability~\cite{ryu-fdtd-quantum-paper-2016,gedney-2011}, we interpret the CFL limit as a condition under which the P-FDTD equations for the region conserve numerical energy. This is desirable, since the conservation of energy is an important property of the continuous equations~\eqref{eq:phi-state-space}--\eqref{eq:a-state-space}. Discretization methods, in general, are not guaranteed to preserve this property~\cite{moerloose-1995-motl}. Additionally, the conservation of numerical energy allows connecting the P-FDTD region to other energy conserving regions without causing numerical instability.
To the best of our knowledge, the connection between the conservation of energy and time step has only been shown for the fields-based FDTD~\cite{edelvik2004general,jnl-2017-fdtd-dissipative,jnl-2018-fdtd-3d-dissipative}. 

\subsection{Vector Potential}

For the vector potential, the storage function and supply rate have the same form \eqref{eq:storage} and \eqref{eq:supply}, but with $\mathbf R$, $\mathbf S$, $\mathbf x$, $\mathbf u$, and $\mathbf y$ resulting from the formulation for the vector potential similar to~\eqref{eq:dyn-sys-a}--\eqref{eq:dyn-sys-b}.
The resulting storage function and supply rate can be seen as discretizations of the following continuous expressions for stored energy~\cite{chew-qem-new-look-2} and supplied energy
\begin{equation}
\label{eq:storage-cont-a}
\mathcal E_{\text{cont}} = \iiint_V \left[
\frac{\mu^{\!-1}\!\!}{2} \vec B\bigcdot \vec B 
+ \frac{\varepsilon}{2} \dfrac{\partial\vec A}{\partial t} \bigcdot \dfrac{\partial\vec A}{\partial t}
+ \frac{\chi}{2} \kappa^2
\right]
dV
\end{equation}
\begin{equation}
	\label{eq:supply-cont-a}
s_{\text{cont}}=\int_{\Delta t}
\oiint_S \left(\mu^{\!-1} \vec B \times \dfrac{\partial \vec A}{\partial t} + \varepsilon \kappa \dfrac{\partial \vec A}{\partial t}\right) \bigcdot (-\hat n)
dS \ dt
\end{equation}
which obey the conservation law in the continuous case~\cite{backhaus-schafer-energy}.
Similarly to the scalar potential case, the CFL limit~\eqref{eq:cfl} is a sufficient condition for the conservation of energy of the dynamical system corresponding to the vector potential.

\section{Numerical Example}
\label{sec:numerical-example}

\begin{figure}[t!]
	\centering
	\includegraphics[]{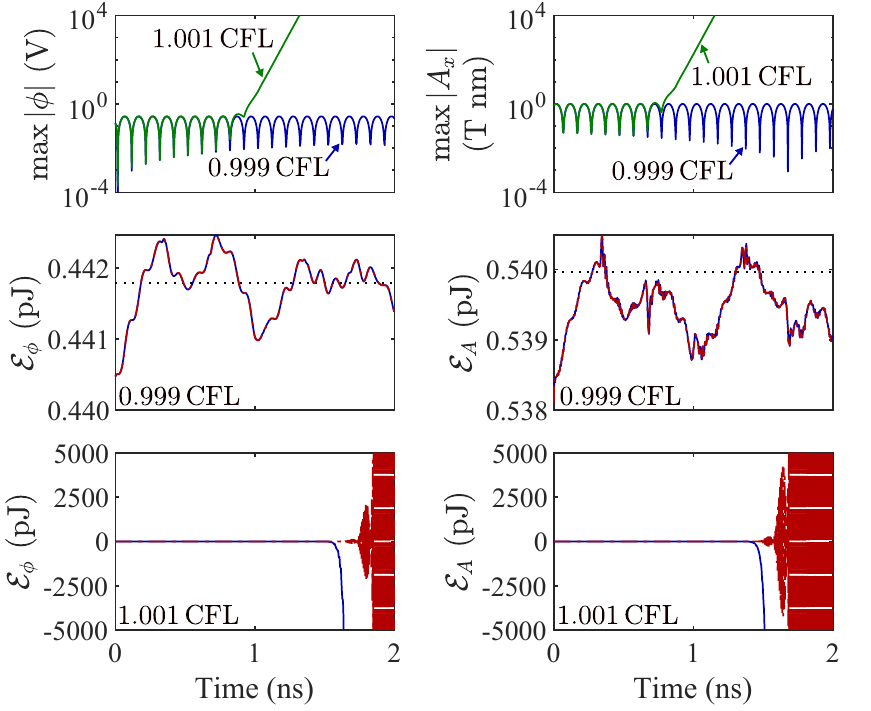}
	\caption{Results of the tests in Sec.~\ref{sec:numerical-example}. Top row: maximum absolute values of potentials over the region. Middle row: storage function~\eqref{eq:storage} ($\color{darkBlue}\st{~~~~}$) and the sum of initial and supplied energy ($\color{darkRed}\st{~~}~\st{~}$) for $\dt=0.999 \dt_{\text{CFL}}$. The dotted black lines show the exact values of energy stored in the scalar and vector potential, evaluated from \eqref{eq:cavity-energy}.
	Last row: storage function~\eqref{eq:storage} ($\color{darkBlue}\st{~~~~}$) and the sum of initial and supplied energy ($\color{darkRed}\st{~~}~\st{~}$) for $\dt=1.001 \dt_{\text{CFL}}$.
	}
	\label{fig:results}
\end{figure}

We verify the theoretical results in Sec.~\ref{sec:dissipativity} by considering a cubic region with side length $a = 10~\text{cm}$ filled with air. We consider a particular solution of the wave equation in the region given by
\begin{equation}
\label{eq:solution-cavity-phi}
\phi = C_{\phi}  \sin(k_x x) \sin(k_y y) \sin(k_z z) \cos(\omega t+\tfrac{\pi}3)
\end{equation}
\begin{equation}
\label{eq:solution-cavity-a}
\vec A = C_A \cos(k_x x) \sin(k_y y)  \sin(k_z z) \sin(\omega t+\tfrac{\pi}3) \ \hat x
\end{equation}
\begin{equation}
k_x = \frac{3 \pi}{a} ,\,
k_y = k_z = \frac{\pi}{a} ,\,
\omega = \sqrt{\tfrac{\varepsilon_0}{\chi_0}(k_x^2+k_y^2+k_z^2)}
\end{equation}
where $C_{A} = 10^{-9}~\text{Tm}$, $C_{\phi} = -C_A \varepsilon_0 k_x/(\chi_0 \omega)$,  $\varepsilon_0$ and $\chi_0$ are values of free-space $\varepsilon$ and $\chi$. 
From \eqref{eq:storage-cont-phi}, \eqref{eq:storage-cont-a}, \eqref{eq:solution-cavity-phi}, and \eqref{eq:solution-cavity-a}, one can calculate the energy stored in the region for the scalar and vector potential
\begin{equation}
	\label{eq:cavity-energy}
	\mathcal E_{\phi,\text{cont}} = \tfrac1{16}\chi_0 C_{\phi}^2 \omega^2 a^3,
	\,\,
	\mathcal E_{A,\text{cont}} = \tfrac1{16}\varepsilon_0 C_A^2 \omega^2 a^3\,.
\end{equation}
The region was divided into $n_x\!=\!90$ and $n_y \!=\! n_z \!=\! 30$ cells. The potentials were updated based on the equations in Sec.~\ref{sec:state-space-form}. The initial conditions on the variables were set in accordance with~\eqref{eq:solution-cavity-phi}--\eqref{eq:solution-cavity-a}. 
The hanging variables on the boundary ($[\nabla \phi]_{\bot}$, $\vec B_{\parallel}$, and $\vec A_{\bot}$) were set by solving for their values from~\eqref{eq:solution-cavity-phi}--\eqref{eq:solution-cavity-a}. For instance, the $[\partial_x\phi]$ on the $i=1$ boundary was found by differentiating~\eqref{eq:solution-cavity-phi} with respect to $x$ and evaluating the result at $x=0$.

The top row of Fig.~\ref{fig:results} shows the evolution of the potentials over time for time steps 0.1\% below and 0.1\% above the CFL limit~\eqref{eq:cfl}. The simulation run below the CFL limit shows no sign of instability. The simulation performed above the CFL limit is unstable, as a result of the algorithm violating the principle of energy conservation.
The two plots in the middle row of Fig.~\ref{fig:results} confirm that the storage functions for scalar and vector potentials are approximately equal to their continuous counterparts \eqref{eq:cavity-energy}.

From \eqref{eq:dissip-ineq}, the storage function should be equal to the energy initially present in the region at the beginning of the simulation plus the total energy absorbed through the boundaries. This is consistent with observation for the test performed under the CFL limit in the second row of Fig.~\ref{fig:results}, where the two curves match with a maximum error comparable to machine precision: $1.8\times 10^{-15}$~pJ and $2.1\times 10^{-15}$~pJ for the scalar and vector potentials, respectively. Interestingly, when the CFL limit is violated (the last row of Fig.~\ref{fig:results}), the two curves deviate substantially, showing that the numerical algorithm violates energy conservation.

\section{Conclusion}

We showed that the update equations of the P-FDTD method~\cite{ryu-fdtd-quantum-paper-2016} can be cast in the form of a discrete-time dynamical system. We proposed suitable discrete expressions to compute the energy stored in the scalar and vector potential present in a region, as well as the power absorbed by the region through its boundary. Through these equations, we identified the precise condition under which the update equations of this method conserve energy in discrete time domain, and discussed its relation to the CFL stability limit.

The proposed developments have multiple uses. The proposed expressions for energy and absorbed power can be used to reliably estimate these quantities, which are useful in characterizing resonators, antennas, and superconductive devices. This study also provided insight into how electromagnetic energy is approximated and behaves in FDTD-like methods based on potentials. Finally, conditions for energy conservation can be used to devise new FDTD schemes based on potentials with guaranteed stability.



%

%

\ifCLASSOPTIONcaptionsoff
  \newpage
\fi



\bibliographystyle{IEEEtran}
\bibliography{IEEEabrv,bibliography}
%

%








\end{document}